\documentclass[12pt]{amsart}
\usepackage{amssymb}

\begin{document}
\begin{center}
\Large Four Entries for\\
Kluwer Encyclopaedia of Mathematics\\
\large Frank Sottile\\
2 June 2000
\end{center}

{\bf The Schubert Calculus}
is a formal calculus of symbols representing
geometric conditions used to solve problems in enumerative geometry.
This originated in work of Chasles~\cite{Ch1864} on conics and
was systematized and used to great effect by Schubert in his treatise ``Kalk\"ul
der abz\"ahlenden Geometrie''~\cite{Sch1879}.
The justification of Schubert's enumerative calculus and the verification of the
numbers he obtained was the 15th problem of Hilbert.

Justifying Schubert's enumerative calculus was a major theme of 20th century
algebraic geometry and {\bf Intersection Theory} provides a satisfactory modern
framework.
{\bf Enumerative Geometry} deals with the second part of Hilbert's problem.
Fulton's book~\cite{Fu98} is a complete reference for Intersection Theory;
for historical surveys and a discussion of Enumerative Geometry, see
the surveys~\cite{Kl76,Kl87}.

The Schubert calculus also refers to 
mathematics arising from the following class of enumerative geometric
problems:
Determine the number of linear subspaces of projective space that satisfy
incidence conditions imposed by other linear subspaces.
For a survey, see~\cite{MR48:2152}.
For example, how many lines in projective 3-space 
meet 4 given lines?
These problems are solved by studying both the geometry and the cohomology
or Chow rings of Grassmann varieties.
This field of Schubert calculus enjoys important connections not only to algebraic
geometry and algebraic topology, but also to algebraic combinatorics,
representation theory, differential geometry, linear algebraic groups, and
symbolic computation, and has found applications in 
numerical homotopy continuation~\cite{HSS98}, linear algebra~\cite{Fu00} and
systems theory~\cite{By80}.
 
The Grassmannian $G_{m,n}$ of $m$-dimensional subspaces ($m$-planes) in 
${\mathbb P}^n$ over a field $k$ has distinguished Schubert varieties
$$
  \Omega_{a_0,a_1,\ldots,a_m} V_{\bullet}\ :=\ 
  \{W\in G_{m,n}\mid W\cap V_{a_j}\geq j\}\,,
$$
where $V_{\bullet}:V_0\subset V_1\subset\cdots\subset V_n={\mathbb P}^n$ is a 
{\bf flag} of linear subspaces with $\dim V_j=j$.
The {\bf Schubert cycle} $\sigma_{a_0,a_1,\ldots,a_n}$ is the cohomology class
Poincar\'e dual to the fundamental homology cycle of
$\Omega_{a_0,a_1,\ldots,a_m} V_{\bullet}$.
The {\bf Basis Theorem} asserts that the Schubert cycles  
form a basis of the Chow ring $A^*G_{m,n}$ (when $k$ is the complex numbers, the
integral cohomology groups $H^*G_{m,n}$) of the Grassmannian
with $\sigma_{a_0,a_1,\ldots,a_m}\in
  A^{(m+1)(n+1)-\binom{m+1}{n+1}-a_0-\cdots-a_m}G_{m,n}$.
The {\bf Duality Theorem} asserts that the basis of
Schubert cycles is self-dual under the {\bf intersection pairing}
$$
   (\alpha,\beta)\in H^*G_{m,n}\otimes H^*G_{m,n}\ \longmapsto\ 
   \deg(\alpha\cdot\beta)=\int_{G_{m,n}} \alpha\cdot \beta\,.
$$
with $\sigma_{a_0,a_1,\ldots,a_m}$ dual to $\sigma_{n-a_m,\ldots,n-a_1,n-a_0}$.

Let $\tau_b:=\sigma_{n-m-b,n-m+1,\ldots,n}$, a {\bf special Schubert cycle}.
Then
$$
  \sigma_{a_0,a_1,\ldots,a_m}\cdot \tau_b\ =\ 
   \sum \sigma_{c_0,c_1,\ldots,c_m}\,,
$$
the sum over all $(c_0,c_1,\ldots,c_m)$ with 
$0\leq c_0\leq a_0\leq c_1\le a_1\leq \cdots\leq c_m\leq a_m$
and $b= \sum_i(a_i-c_i)$.
This {\bf Pieri Formula} determines the ring structure of cohomology;
an algebraic consequence is the {\bf Giambelli formula} for expressing an arbitrary
Schubert cycle in terms of special Schubert cycles. 
Define $\tau_b=0$ if $b<0$ or $b>m$, and $\tau_0=1$.
Then Giambelli's formula is
$$
  \sigma_{a_0,a_1,\ldots,a_m}\ =\ 
    \det [ \tau_{n-m+j-a_i}]_{i,j=0,\ldots,m}\,.
$$

These four results enable computation in the Chow ring
of the Grassmannian, and the solution of many problems in enumerative geometry.
For instance, the number of $m$-planes meeting $(m+1)(n-m)$ general
$(n-m-1)$-planes
non-trivially is the coefficient of $\sigma_{0,1,\ldots,m}$ in the product
$(\tau_1)^{(m+1)(n-m)}$, is~\cite{Sch1886c}
$$
    \frac{1!\,2!\cdots(n\!-\!m\! -\!1)!\cdot[(m+1)(n-m)]!}
  {(n\!-\!m)!\,(n\!-\!m\! +\!1)!\cdots (n! -1)!}\ . 
$$  

These four results hold more generally for cohomology rings of flag manifolds $G/P$; 
Schubert cycles form a self-dual basis, the Chevalley formula~\cite{Ch91}
determines the ring structure (when $P$ is a Borel subgroup), and
the formulas of Bernstein-Gelfand-Gelfand~\cite{BGG73} and Demazure~\cite{De74}
give the analog of the Giambelli formula.
More explicit Giambelli formulas are provided by {\bf Schubert polynomials}.

One cornerstone of the Schubert calculus for the Grassmannian
is the {\bf Littlewood Richardson rule}~\cite{LR34} for expressing a
product of Schubert cycles in terms of the basis of Schubert cycles.
[This rule is usually expressed in terms of an alternative indexing of Schubert
cycles using partitions.
A sequence $(a_0,a_1,\ldots,a_m)$ corresponds to the partition
$(n-m-a_0,n-m+1-a_1,\ldots,n-a_m)$.]
The analog of the Littlewood Richardson rule is not known for most other flag
varieties $G/P$.\medskip

\bigskip

\hrule
\hrule
\bigskip

{\bf Schubert Cell} --
The orbit of a Borel subgroup $B\subset G$ on a flag variety
$G/P$~\cite[14.12]{Bo91}.
Here, $G$ is a semisimple linear algebraic group over an algebraically closed
field $k$ and $P$ is a parabolic subgroup of $G$ so that $G/P$ is a complete
homogeneous variety.
Schubert cells are indexed by the cosets of the Weyl group $W_P$ of $P$ in the
Weyl group $W$ of $G$.
Choosing $B\subset P$, these cosets are identified with $T$-fixed points of $G/P$,
where $T$ is a maximal torus of $G$ and $T\subset B$.
The fixed points are conjugates $P'$ of $P$ containing $T$.
The orbit $BwW_P\simeq {\mathbb A}^{\ell(wW_P)}$, the affine space of dimension
equal to the length of the shortest element of the coset $wW_P$.
When $k$ is the complex numbers, Schubert cells constitute a CW-decomposition
of $G/P$.

Let $k$ be any field and suppose $G/P$ is the Grassmannian $G_{m,n}$ of $m$-planes in
$k^n$. 
Schubert cells for $G_{m,n}$ arise in an elementary manner.
Among the $m$ by $n$ matrices whose row space is a given $H\in G_{m,n}$, there is
a unique echelon matrix
$$
 \left[
  \begin{array}{ccccccccccccccc}
   *&\cdots&*&1&0&\cdots&0&0&0&\cdots&0&0&0&\cdots&0\\
   *&\cdots&*&0&*&\cdots&*&1&0&\cdots&0&0&0&\cdots&0\\
   \vdots&\ddots&\vdots&\vdots&\vdots&\ddots&\vdots&
   \vdots&\vdots&\ddots&\vdots&\vdots&\vdots&\ddots&\vdots\\
   *&\cdots&*&0&*&\cdots&*&0&*&\cdots&*&1&0&\cdots&0
  \end{array}
 \right]
$$
This echelon representative of $H$ is computed from any representative by Gaussian
elimination. 
The column numbers $a_1<a_2<\cdots<a_m$ of the leading entries (1s) of the rows
in this echelon representative determines the type of $H$.
Counting the undetermined entries in such an echelon matrix shows that 
the set of $H\in G_{m,n}$ with this type is isomorphic to 
${\mathbb A}^{mn-\sum(a_i+i-1)}$.
This set is a Schubert cell of $G_{m,n}$.\medskip

\bigskip

\hrule
\hrule
\bigskip

{\bf Schubert Cycle} --
The cycle class of a Schubert variety in the cohomology ring of a complex
flag manifold $G/P$, also called a {\bf Schubert class}.
Here, $G$ is a semisimple linear algebraic group and $P$ is a parabolic subgroup. 
Schubert cycles form a basis for the cohomology
groups~\cite{Sch1886b}~\cite[14.12]{Bo91} of $G/P$. 
They arose~\cite{Sch1886b} as representatives of Schubert conditions on linear
subspaces of a vector space in Schubert's calculus for 
enumerative geometry~\cite{Sch1879}.
The justification of Schubert's calculus in this context by
Ehresmann~\cite{Eh34} realized Schubert cycles as cohomology classes Poincar\'e
dual to the fundamental homology cycles of Schubert varieties in the
Grassmannian. 
While Schubert, Ehresmann, and others worked primarily on the Grassmannian,
the pertinent features of the Grassmannian extend to general flag varieties
$G/P$, giving Schubert cycles as above.

More generally, when $G$ is a semisimple linear algebraic group over a field, 
there are Schubert cycles associated to Schubert varieties in each of the
following theories for $G/P$:
singular (or deRham) cohomology, the Chow ring, K-theory, or
equivariant or quantum versions of these theories.
For each, the Schubert cycles form a basis over the base ring.
For the cohomology or the Chow ring, the Schubert cycles are universal characteristic
classes for (flagged) $G$-bundles.
In particular, certain {\bf special Schubert cycles} for the Grassmannian are
universal Chern classes for vector bundles.\medskip

\bigskip

\hrule
\hrule
\bigskip

{\bf Schubert Polynomials}
were introduced by Lascoux and Sch\"utzenberger~\cite{LS82a} as distinguished
polynomial representatives of Schubert cycles in the cohomology ring of the
manifold $F\ell_n$ of complete flags in ${\mathbb C}^n$.
This extended work of Bernstein-Gelfand-Gelfand~\cite{BGG73} and
Demazure~\cite{De74} who gave algorithms for computing representatives of Schubert
cycles in the coinvariant algebra, which is isomorphic to the cohomology ring of
$F\ell_n$~\cite{Bo53}
$$
  H^*(F\ell_n,{\mathbb Z})\ \simeq\ {\mathbb Z}[x_1,x_2,\ldots,x_n]/
      {\mathbb Z}^+[x_1,x_2,\ldots,x_n]^{{\mathcal S}_n}\,.
$$
Here, ${\mathbb Z}^+[x_1,x_2,\ldots,x_n]^{{\mathcal S}_n}$ is the ideal generated
by the non-constant polynomials that are symmetric in $x_1,x_2,\ldots,x_n$.
Macdonald~\cite{Mac91} gives an elegant algebraic treatment of Schubert
polynomials, while Fulton~\cite{Fu97} and Manivel~\cite{Man98} deal more with
geometry.

For each $i=1,2,\ldots,n-1$, let $s_i$ be the transposition $(i,i+1)$ in the
symmetric group ${\mathcal S}_n$, which acts on 
${\mathbb Z}[x_1,x_2,\ldots,x_n]$.
The {\bf divided difference} operator $\partial_i$ is defined by
$$
  \partial_i f\ =\ (f- s_i f)/(x_i-x_{i+1})\,.
$$
These satisfy 
\begin{equation}\label{eq:dd}
  \begin{array}{rcl}
   \partial_i^2&=&0\\
   \partial_i \partial_j&=&\partial_j \partial_i\qquad\qquad\ 
       \mbox{ if }\ |i-j|>1\\
   \partial_i\partial_{i+1}\partial_i&=&\partial_{i+1}\partial_i\partial_{i+1}
  \end{array}
\end{equation}
If $f_w\in{\mathbb Z}[x_1,x_2,\ldots,x_n]$ is a representative of the Schubert
cycle $\sigma_w$, then 
$$
  \partial_i f_w\ =\ \left\{\begin{array}{lcl}
           0&\quad&\mbox{if }\ell(s_iw)>\ell(w)\\
       f_{s_iw}&&\mbox{if }\ell(s_iw)<\ell(w)
      \end{array}\right.\ ,
$$
where $\ell(w)$ is the length of a permutation $w$ and $f_{s_iw}$ represents
the Schubert cycle $\sigma_{s_iw}$.
Given a fixed polynomial representative of the Schubert cycle $\sigma_{w_n}$
(the class of a point as $w_n\in{\mathcal S}_n$ is the longest element),
successively applying divided difference operators gives polynomial
representatives of all Schubert cycles, which are independent of the choices
involved, by~(\ref{eq:dd}).

The choice of representative 
${\mathfrak S}_{w_n}=x_1^{n-1}x_2^{n-2}\cdots x_{n-1}$ for $\sigma_{w_n}$
gives the {\bf Schubert polynomials}.
Since $\partial_n\cdots\partial_1{\mathfrak S}_{w_{n+1}}={\mathfrak S}_{w_n}$,
Schubert polynomials are independent of $n$ and give polynomials
${\mathfrak S}_w\in{\mathbb Z}[x_1,x_2,\ldots]$ for 
$w\in{\mathcal S}_\infty=\bigcup {\mathcal S}_n$.
These form a basis for this polynomial ring, and every {\bf Schur polynomial}  
is also a Schubert polynomial.

The {\bf transition formula} gives another recursive construction of Schubert
polynomials.
For $w\in{\mathcal S}_\infty$, let $r$ be the last descent of $w$ 
($w(r)>w(r+1)<w(r+2)<\cdots$) and define $s>r$ by $w(s)<w(r)<w(s+1)$.
Set $v=w(r,s)$, where $(r,s)$ is the transposition.
Then
$$
  {\mathfrak S}_w\ =\  x_r {\mathfrak S}_v + \sum {\mathfrak S}_{v(q,r)}\,,
$$
the sum over all $q<r$ with $\ell(v(q,r))=\ell(v)+1=\ell(w)$.
This formula gives an algorithm to compute ${\mathfrak S}_w$
as the permutations that appear on the right hand side are either shorter than 
$w$ or precede it in reverse lexicographic order, and the minimal such permutation
$u$ of length $m$ has ${\mathfrak S}_u=x_1^m$.

The transition formula shows that the Schubert polynomial ${\mathfrak S}_w$
is a sum of monomials with {\it nonnegative} integral coefficients.
There are several explicit formulas for the coefficient of a monomial in a Schubert
polynomial, either in terms of the weak order of the symmetric
group~\cite{BJS93,Be92,FS95}, an intersection number~\cite{KM97}, or the Bruhat
order~\cite{BS00}. 
An elegant conjectural formula of Kohnert~\cite{Kohnert} remains unproven.
The Schubert polynomial ${\mathfrak S}_w$ for $w\in{\mathcal S}_n$ is also the
normal form reduction of any polynomial representative of the Schubert cycle
$\sigma_w$  with respect to the degree reverse lexicographic 
term order on ${\mathbb Z}[x_1,x_2,\ldots,x_n]$ with $x_1<x_2<\cdots<x_n$.

The above-mentioned results of~\cite{Bo53,BGG73,De74} are valid more generally for 
for any flag manifold
$G/B$ with $G$ a semisimple reductive group and $B$ a Borel subgroup.
When $G$ is an orthogonal or symplectic group, there are competing theories of
Schubert polynomials~\cite{BH95,FK96,LPR98}, each with own merits.
There are also double Schubert polynomials suited for computations of
degeneracy loci~\cite{FP98}, quantum Schubert polynomials~\cite{FGP97,C-F99},
and universal Schubert polynomials~\cite{Fu99}.\medskip


\providecommand{\bysame}{\leavevmode\hbox to3em{\hrulefill}\thinspace}

\hfill {\it Frank Sottile}

\hfill {\it email address}: {\tt sottile@math.umass.edu}

\end{document}